\documentclass{article}
\usepackage{arxiv}
\usepackage[utf8]{inputenc} % allow utf-8 input
\usepackage[T1]{fontenc}    % use 8-bit T1 fonts
\usepackage{hyperref}       % hyperlinks
\usepackage{url}            % simple URL typesetting
\usepackage{booktabs}       % professional-quality tables
\usepackage{amsfonts}       % blackboard math symbols
\usepackage{nicefrac}       % compact symbols for 1/2, etc.
\usepackage{microtype}      % microtypography
\usepackage{graphicx}
\usepackage{float}
\usepackage{doi}
\usepackage[linesnumbered,ruled,vlined]{algorithm2e}
\usepackage{algpseudocode}
\usepackage{subcaption}
\usepackage{amsmath}
\usepackage{amsfonts}
\usepackage{amssymb}
\usepackage{amsthm}
\usepackage{lineno}
\usepackage{tabularx}
\usepackage{array}
\usepackage[english]{babel}

\newtheorem{definition}{Definition}

\newtheorem{remark}{Remark}
\DeclareMathSymbol{\N}{\mathbin}{AMSb}{"4E}
\DeclareMathSymbol{\Z}{\mathbin}{AMSb}{"5A}
\DeclareMathSymbol{\R}{\mathbin}{AMSb}{"52}
\DeclareMathSymbol{\Q}{\mathbin}{AMSb}{"51}
\DeclareMathSymbol{\C}{\mathbin}{AMSb}{"43}
\DeclareMathSymbol{\T}{\mathbin}{AMSb}{"54}

\title{Estimating Varying Parameters in Dynamical Systems: A Modular Framework Using Switch Detection, Optimization, and Sparse Regression}

\date{January 2, 2025}

\author{ 
	\href{https://orcid.org/0000-0002-6535-1579}           {\includegraphics[scale=0.06]                          {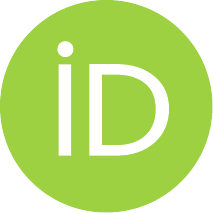}\hspace{1mm}Jamiree Harrison*} \\
	Department of Mechanical Engineering\\
	University of California, Santa Barbara\\
	Santa Barbara, California \\
	\texttt{jamiree@ucsb.edu} \\
	\And
	\href{https://orcid.org/0000-0001-7630-7429}           {\includegraphics[scale=0.06]                          {orcid.pdf}\hspace{1mm}Enoch Yeung} \\
	Department of Mechanical Engineering\\
	University of California, Santa Barbara\\
	Santa Barbara, California \\
	\texttt{eyeung@ucsb.edu} \\
}

\renewcommand{\shorttitle}{\textit{arXiv}}

\begin{document}
\maketitle

\begin{abstract}
The estimation of static parameters in dynamical systems is well-established, but challenges remain for systems with parameters that vary over time or space. Assuming prior knowledge of the system’s dynamics, we aim to identify functions describing these varying parameters. Initially, we focus on cases where parameters are discretely switching piece-wise constant functions. For this sub-class, we propose an algorithmic framework to detect discrete parameter switches and fit a piece-wise constant model using optimization-based parameter estimation. 
Our framework is modular, allowing customization of switch detection, numerical integration, and optimization steps. We demonstrate its utility through examples, including a time-varying promoter-gene expression model, a genetic toggle switch, parameter-switching manifolds, mixed fixed and switching parameters, non-uniform switching, and models governed by partial differential equations like the heat and advection-diffusion equations. Binary segmentation is used for switch detection, and Nelder-Mead and Powell methods for optimization. 
Additionally, we assess the framework's robustness to measurement noise. For more complex cases, we extend the approach using dictionary-based sparse regression with trigonometric and polynomial functions to capture continuously varying parameters. This comprehensive framework offers a versatile solution for parameter estimation in varying dynamical systems.
\end{abstract}


\section{Introduction}

Parameter estimation is a fundamental inverse problem which finds applications in control theory, statistics, systems biology, and physics, to name just a few. Estimating a system's parameters often precedes the prediction of that system's future states, the understanding of that system's dynamics, and the structure of the given system. Obtaining a system's parameters in certain instances is a key step in the accurate control of that system.  To exemplify the estimation of a
system's parameters for its accurate control, the work of Ryan and Ichikawa \cite{RYAN1979149} shows how this can be done in the context of a linear distributed system with static parameters that get updated estimates. Furthermore, parameter estimation has been shown to provide benefits in addressing model discrepancies in control systems.  For example, Burns et al. \cite{BURNS201411679} address the accurate control of systems with delays, model errors, and measurement errors by employing hierarchical modeling to estimate uncertain disturbances, offering an alternative to Bayesian analysis.

Mitra and Hlavacek offer a comprehensive review of parameterization and uncertainty quantification in mathematical models of immunoreceptor signaling and other biological processes. They highlight techniques such as gradient-based estimation, profile likelihood, bootstrapping, and Bayesian inference, emphasizing their applications in immune system modeling \cite{Mitra_Hlavacek_2019}. Penas et al. introduced saCeSS, a novel parallel optimization method for parameter estimation in large-scale kinetic models. This approach demonstrated robust performance, significantly reduced computation times, and user-friendly self-tuning, with potential applications in systems biology and whole-cell dynamic modeling \cite{Penas}.

Zhan and Yeung proposed two innovative parameter estimation methods that combine spline theory with Linear and Nonlinear Programming. These methods offer robust, efficient, and flexible solutions for identifying parameters in systems biology models without the need for ODE solvers \cite{Zhan_Yeung_2011}. Valderrama-Bahamóndez and Fröhlich compared Markov Chain Monte Carlo (MCMC) techniques for estimating rate parameters in non-linear ODE systems used in systems biology. They found that parallel adaptive MCMC, combined with informative Bayesian priors, outperformed other methods for high-dimensional, multi-modal parameter distributions \cite{Valderrama}. Giampiccolo et al. introduced a hybrid Neural ODE framework for parameter estimation and identifiability assessment in biological models with partial structural knowledge. Their approach integrates hyperparameter tuning with \textit{a posteriori} identifiability analysis to address challenges posed by neural networks, demonstrating effectiveness on noisy, real-world-inspired benchmarks \cite{Giampiccolo_Reali_Fochesato_Iacca_Marchetti_2024}.

For additional examples in quantitative biology, Ashyraliyev et al. provide a mini-review exploring the applications of mathematical models in biology. They focus on parameter inference, model identifiability, and commonly used methods for parameter space exploration \cite{ashraliyev}. Beyond quantitative biology, parameter estimation has been successfully applied in fields like astrophysics. For instance, Douspis et al. \cite{Douspis2001} and Bailer-Jones \cite{bailerjones} have utilized parameter estimation across various models and datasets, demonstrating its versatility and effectiveness.

The problem we aim to solve in this paper is at the interface dynamical systems theory and parameter estimation. Dynamical systems in which the parameters vary as the system evolves are known as parameter-varying systems, and these parameters can change continuously or switch discretely.  Systems that are modeled with the prior assumption of having static parameters could often be better characterized, better controlled, and more accurately captured via parameter-varying models, and our paper serves as a novel avenue for fitting these more accurate models.  An equivalent classification of this set of problems is the input estimation problem because mathematically representing these time-varying parameters as inputs which change over time is indistinguishable as noted in \cite{martinelli2024}. 
Accurately estimating the varying parameters along with their swtich times for general classes of systems remains a nontrivial task. In this paper, we address this challenge.

\subsection{Definitions}
We first need to define what constitutes a discrete parameter switch in a parameter-varying system, as well as clarify what it means for a parameter to vary continuously. 

\begin{definition}\label{def: cts param} \textbf{Continuously Varying Parameter}

We say that $p(t):\R \rightarrow \R$ is a \textbf{continuously varying parameter} if $p(t)$ is continuous and not a constant function. 
    
\end{definition}

\begin{definition}\label{def: param switch} \textbf{Discrete Parameter Switch}

    Let $n_p$ be defined as the number of parameters in a given system.  
    Suppose that we have piece-wise parameters $p(t): \R \rightarrow \R^{n_p}$, given by 
    \begin{align}
        p(t) = 
        \begin{cases}
             p_1(t)  &  \text{if} \quad t \in [t_0, t_1] \\
             p_2(t)  &  \text{if} \quad t \in (t_1, t_2] \\
                  & \vdots \\
             p_n(t)  &  \text{if } \quad t \in (t_{n-1}, t_n] 
        \end{cases}
    \end{align}
    where $n \in \N$ and $p_i(t)$ are continuously varying .  If $p_k (t_k) \neq p_{k+1} (t_k)$ for a given $k \in \{ 1,...,n-1 \}$, then we say that a \textbf{discrete parameter switch} has occured at $t_k$.  The number of discrete parameter switches is the \textbf{discrete switch number} $N_{s_d}$.  These \textbf{discrete switch locations} which we will denote by $t_{k_d}$ are jump discontinuities of $p(t)$.
    
    We say that a \textbf{continuous parameter switch} occurs at \textbf{continuous switch locations} denoted by $t_{k_c}$ at values of $t$ where $p(t)$ is continuous but not differentiable.  Explicitly, a continuous switch location occurs if $p_k (t) \neq p_{k+1} (t)$ for some $t \in [t_{k-1}, t_{k+1}] \setminus t_k$, but $p_k (t_k) = p_{k+1} (t_k)$.  The number of continuous parameter switches is the \textbf{continuous switch number} $N_{s_c}$.  The \textbf{total switch number} is defined as $N_s:= N_{s_c} + N_{s_d}$ . 
\end{definition}

\begin{remark}
    Suppose that we have a dynamical system $\frac{dX}{dt}=f(X,p(t))$, and suppose that $f$ is continuous with respect to $t$ if $p(t)$ is continuous (without discrete switches).  
    
    If we have $\frac{dX}{dt}=f(X,p(t))$ such that $p(t)$ represents piece-wise parameters with discrete parameter switches at $t_k$ (and thus, jump discontinuities) as described in Definition (\ref{def: param switch}), then $f$ has jump discontinuities at $t_k$ and $X(t)$ is non-differentiable at $t_k$.   
\end{remark}
The above remark expresses that the discrete parameter switches of a parameter-varying system manifest in a given trajectory, $X(t)$, as non-differentiable points in $t$.  In a data-driven empirical setting where we are given discrete data points, $X_{data}$, we cannot analytically determine where these parameter switch locations are, so we leverage the data-driven switch detection methods that are available in the literature (see Section \ref{sec: binseg}) and we extend their implementations to detect switching parameters of dynamical systems.

\subsection{General Problem Formulation (continuous-time systems)}
\label{sec: prob form}
Assume that we have data denoted $X_{data} \in \R^{m \times N}$ from a parameter-varying system denoted 
\begin{align}
    \frac{dX}{dt}=f(X,p(t))
\end{align}
with $m \in \N$ states, $N \in \N$ data points, vector field $f: \R^m \rightarrow \R^m$, and $n_p \in \N$ parameters $p(t): \R \rightarrow \R^{n_p}$. 
    
Further, assume we know the form of the system's model, $f(X,p(t))$ but we do not know the varying parameters $p(t)$.  Lastly, assume that the parameters, $p(t)$, are piece-wise functions.
    
\textbf{Problem}: Identify $p(t)$ given $X_{data}$,  i.e. find the minimizing parameters such that 
    \begin{align}
        p(t) = \textbf{argmin} \|X_{data} - X_{model}(p(t)) \|_2
    \end{align}
    where  $X_{model} \in \R^{m \times N}$ 
 is given by  $X_{model} 
 = \int\limits_{\text{numerical}}  \frac{dX}{dt} dt$ .

\begin{remark}
    Here, the specified numerical integration scheme (ODE solver) is presupposed in Assumption 2 such that it will grant sufficiently accurate $X_{model}$ data in terms of $\|X_{model} - X(t_{d}) \|_2 < \epsilon$ for some desired $\epsilon>0$ where $X(t_d)$ is the corresponding discretization of $X(t)$, the solution to the system. This problem formulation also accounts for partial differential equations (PDEs) and discrete-time parameter-varying systems are accounted for in the supplement. 
\end{remark}

This problem formulation encompasses several sub-problems, each addressed with varying degrees of success in the literature. While these contributions demonstrate rigor and specificity in tackling their respective challenges, our framework addresses the general class of problems holistically, filling critical gaps. Additionally, we acknowledge that the specialized methods discussed in the literature may outperform our framework in certain domains and should be employed when appropriate. However, our approach effectively covers many use cases within the largely unexplored space of the broader problem class.

For example, Wang et al. proposed a method for switch detection and robust identification in slowly switched Hammerstein systems \cite{WANG2019202}. Zheng et al. developed a technique for identifying fast-switched linear-nonlinear alternating subsystems, particularly in high-speed train models \cite{ZHENG2024105815}. Yu et al. \cite{YU2023101364} and Bianchi et al. \cite{BIANCHI2021109415} introduced randomized methods to address switched ARX and NARX systems, SISO systems, and other non-autonomous systems. Li et al. created methods for identifying nonlinear Markov jump systems \cite{LI2022348}.  Boddupalli et al. have developed machine learning methods involving symbolic regression to give interpretable representations of dynamical systems without the need of instantiating known models or dictionary functions, but these methods struggle to identify systems with varying parameters \cite{nibodh, nibodh2}.  These methods stand in contrast to the recent data-driven system identification and control developments utilizing Koopman operator theory frameworks \cite{johnsonnonlinear, susuki2022control, susuki2024control}. 

Kon et al. \cite{kon2023direct} presented an input-output linear parameter-varying (LPV) feedforward parameterization technique and a corresponding data-driven estimation method, leveraging neural networks to model the dependency of coefficients on scheduling signals. Xu et al. developed an algorithm for parameter estimation in time-varying systems with invariant matrix structures \cite{Xu2023}. Similarly, Luiz and Reinaldo et al. tackled parameter estimation in linear, continuous, time-varying dynamical systems, assuming known coefficient matrices, measurable states, and bounded piecewise continuous parameters \cite{CLAUDIOANDRADESOUZA2011777}. 

Johnson et al. \cite{hybrid_param} applied hybrid gradient descent to estimate unknown constant parameters in hybrid systems with flow and jump dynamics affine in the parameters. Qian et al. used Bayesian parameter estimation for discrete-time linear parameter-varying finite impulse response systems \cite{qian}. Altaie et al. \cite{altaie2024} demonstrated how parameter estimation in parameter-varying systems can inform the design of robust controllers for multidimensional systems.  

While accurately detecting a system's changing parameters has been achieved for certain classes of systems most notably with methods developed my Yong et al. and Curi et al. respectively \cite{yong2017simultaneous, curi2019adaptive}, our paper presents a generalized framework for estimating the switching parameters of systems, their switch locations, and for continuously varying parameters. We do not do a comparison study between our framework and the methods of others because there are several permutations of our modular framework that could give different results, and no other method can estimate the switch locations and the parameters as we have done. At a glance, we are able to get similar results as several of the aforementioned methods. 
Further, we showcase implementations of algorithms within our framework that successfully handle several examples including a time-varying promoter-gene expression model, the time-varying genetic toggle switch, a parameter-switching manifold, a model with a mixture of fixed and switching parameters, a model in which parameters are not switching uniformly, and the advection-diffusion equation with a time-varying coefficient of advection.  Our paper contributes to the body of literature in that it provides a framework to estimate the changing parameters of parameter varying dynamical systems using an approach that utilizes the combination of switch detection, optimization, and sparse dictionary regression. 

\subsection{Contributions}

\begin{enumerate}
	\item 
	Ensemble methodological framework that allows for the estimationn of time-varying parameters/inputs for the broad class of parameter/input estimation problems.
	
	\item 
	Precise switch locations of parameter/input estimates obtained from data, allowing for accurate parameter/input based data segmentation.
	
	\item 
	Function parameterizations of time-varying parameter/input estimates are obtained from data with sparse regression.  
	
	\item 
	Adapted signal change-point detection algorithms to detect parameter/input switches in dynamical systems.
\end{enumerate}

\section{Methodology for algorithmic framework}

\begin{figure}[H]
    \centering
    \includegraphics[width=.9\linewidth]{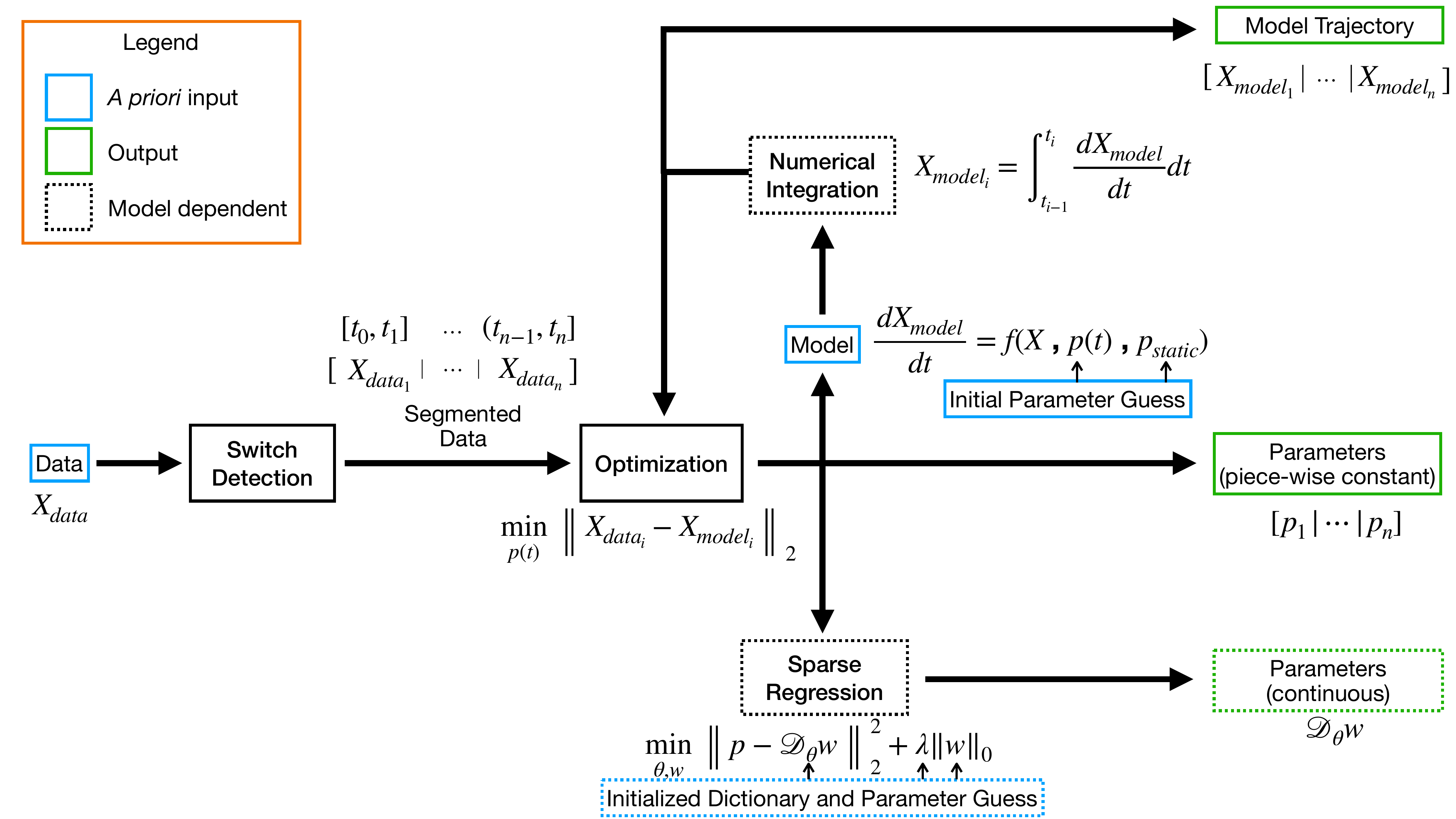}
    \caption{Outline of algorithmic framework.  The data (left) get segmented by a switch detection algorithm.  The presupposed model (middle-right) is then fit to each data interval via minimization of the model trajectory error (middle-left).  If the model is a differential equation, then the numerical integration sub-step (upper-middle) must take place in the optimization loop (middle), but numerical integration need not take place if the model is a parametric equation in the form of $X_{model} = f(X,p(t))$. 
 Optimal parameters and model trajectories are then given as output (right). 
 Sparse regression is used if the parameters of the model are pre-supposed to be continuously varying as opposed to being discretely switching piece-wise constant parameters. } 
    
    \label{fig:algo flow}
\end{figure}

To give a broad overview of our framework, let us define $X_{data} \in \R^{m \times N}$ as the input data from our system of interest, and $X_{model} \in \R^{m \times N}$ as the predictive model data that approximates $X_{data}$ upon completion of an algorithm within our framework.  These are matrices with $m$ being the number of states and $N$ being the number of data points (or snapshots) such that a row of either $X_{data}$ or $X_{model}$ is a state trajectory. 
The overview of our framework is given in Figure \ref{fig:algo flow} . 

The framework in Figure (\ref{fig:algo flow}) starts with the initialized continuous-time model, then detects switches in the input data and then segments the data into $n$ intervals.  From here, model parameters $p(t)$ and $p_{static}$ are identified to fit the model to the input data over each segmented time interval.  This is done by numerically integrating the segmented dynamic model to obtain $X_{model_i}$ and then using numerical optimization to minimize the distance (norm) between $X_{data_i}$ and $X_{model_i}$ where $i \in \{ 1,...,n \}$.  Numerical integration is used for differential equation models in the form of $\frac{dX}{dt} = f(X, p(t))$ and is not used for parametric curve models in the form of $ X_{model} = f(X, p(t))$. If the parameters were presupposed to be piece-wise constant, then our algorithm outputs them as such.  If the parameters were presupposed to be continuously varying functions, then we assume a high number of switches to have occurred in the data and use dictionary-based sparse regression to fit continuous functions to the parameters.   Finally, the identified model parameters are used to generate the model trajectory, $X_{model}$ that accurately approximates $X_{data}$.  Being that this trajectory is an approximation of the $N$ sequential snapshots given in $X_{data}$, these are also known as $N$-step predictions.

From here it is important to note that this framework works for ODEs and PDEs that are written to have a first order derivative in time.  In the context of numerical methods for ODE solving, this means that the ODE must be in standard form. One can handle PDEs within our framework using the method of lines (discretizing the spatial derivative).  Further, we showcase and provide code in the supplement for our given implementations, but it should not be understated that given the modular structure of the code, that a user can swap out the switch detection, optimization, and numerical integration schemes for preferred methods.  
Further, our framework is easily modulated to the estimate the varying parameters of a discrete time system (found in supplement).

%\begin{proposition}\label{Convergence} \textbf{Convergence of Algorithms in framework.}
%
%    An algorithm specified under our framework will converge if the following sub-steps are convergent: 
%    
%    1. Switch detection / data segmentation 
%    
%    2. Numerical integration of parameterized model 
%    
%    3. Minimization of parameterized model error 
%\end{proposition}
%
%\begin{proof}
%    Note that we have a finite number of processes as sub-steps in our algorithmic framework.  Each sub-step is assumed to converge to either a locally or globally optimal solution depending on the method chosen.  The proof is then trivial in that we have a finite sequence of convergent processes, so the sequence of processes itself must converge to either a locally or globablly optimal solution. 
%\end{proof}

\begin{remark}
For globally optimal solutions to be granted under this framework (over finite domains), all sub-steps listed in the framework must be globally optimal.  If global optimization algorithms are used over a finite domain of possible parameters, then globally optimal parameters are guaranteed over that domain, however, these globally optimal parameters may not be unique because other distinct minima of the specified non-convex objective function, $\|X_{model} - X_{data} \|_2$, may exist. 
\end{remark}

\subsection{Evaluative metrics}

Now that we have established our algorithmic framework, we formalize the metrics to be used for the evaluation of such algorithms. The parameter error $E_p$ is defined by the distance between a system's set of true parameters over the discrete data points, $p(t_d)$, and the approximate set of parameters of that system, $\hat{p}$, normalized by the number of parameters, $n_p$. The trajectory error $E_t$ is defined by the distance between the input data matrix, $X_{data}$, and the data matrix generated by the given model with estimated parameters, $X_{model}$, normalized by the number of data points $N$.  From context, it is clear that the trajectory error $E_t$ is the objective function minimized in our framework.  The switch number error $E_{N_s}$ is the difference between the true number of switches, $N_s$, and those detected, $\hat{N_s}$.  Note that if $E_{N_s}$ is positive, then the number of switches detected is an underestimate and if $E_{N_s}$ is negative, then the number of switches detected is an overestimate.  In this way, the sign of $E_{N_s}$ can be used to quantify if a given algorithm systematically under-reports or over-reports the correct number of switches in data.  To measure the distance from the true switches, $\{t_1 ,..., t_{n-1}\}$, to the detected switch estimates, $\{\hat{t}_1 ,..., \hat{t}_{\hat{n}-1}\}$, we use the Hausdorff metric, $H_s$ which measures the worst error among all of the detected switches and is explicitly given by 
\begin{equation}\label{eq: haus}
    H_s := \max \{ \max\limits_{k} \min\limits_{l} |t_k - \hat{t}_l |,  \max\limits_{k} \min\limits_{l} |\hat{t}_k - t_l |  \}
\end{equation}

\begin{remark}
    For explicit clarity, the following would generally be unknown in a practical scenario before any algorithm is applied: the true parameters $p(t_d)$, the true number of switches $N_s$, and the locations of the switches $t_k$.  The point of our framework and the goal of solving the previously stated problem is to infer what these unknowns are.  
    We provide these metrics so that algorithms suited to our stated problem-class can have their performances bench-marked when applied to datasets and systems for which $p(t_d)$, $N_s$, and $t_k$ are known.  
\end{remark}

\renewcommand{\arraystretch}{2}
\begin{table}[!htb]
\begin{center}
\begin{tabular}{ | c | c | c | } 

  \hline
  
  \textbf{Error Type} & \textbf{Description} & \textbf{Formula} \\
  
  \hline

  Parameter error & gap between true parameters and estimated parameters & $E_p : = \frac{\| p(t_d) - \hat{p} \|_2}{n_p}$ \\ 

  \hline
  
  Trajectory error & gap between input data and predicted trajectory from model & $E_t : = \frac{\| X_{data} - X_{model} \|_2}{N}$ \\ 
  
  \hline
  
  Switch number error & difference between true number of switches and those detected & $E_{N_s} := N_s - \hat{N_s} $  \\ 

  \hline
  
   Switch distance &  the worst error among all of the detected switches & $H_s :$ (Equation \ref{eq: haus})  \\
  \hline
\end{tabular}
\end{center}
\caption{Metrics for evaluation of algorithm performance.}
\label{tab: error types}
\end{table}

We now briefly explain our choice of the switch number error \( E_{N_s} \) and the Hausdorff metric \( H_s \) over other metrics commonly reviewed in the literature \cite{truong}. 
The precision and recall metrics, for instance, declare a change point as "detected" if at least one computed change point is within a "margin" of data points from it \cite{truong}. However, these metrics require the introduction of an additional hyperparameter, "margin," which is undesirable for our use case, as it already involves a significant number of user-defined parameters. Another alternative is the Rand index, which measures the similarity between two segmentations as the proportion of agreement between their partitions \cite{randindex}. While this metric has merit, we chose the Hausdorff metric over the Rand index because we find the distance between true and detected switches to be a more intuitive and interpretable measure than a similarity score between the sets of true and detected switches. 
Ultimately, we recognize that the choice of metrics is inherently subjective. Nevertheless, we assert that \( E_{N_s} \), in combination with \( H_s \), provides an effective and practical means of evaluating errors in our specific context.

\subsection{Detection of parameter switches in data}
\label{sec: binseg}

Switch detection (commonly referred to as change point detection) in noisy signals remains challenging, and the optimal method for a given dataset depends on the signal type, noise level, and number of expected changes. One popular method, Bayesian Online Change Point Detection (BOCPD), employs a Bayesian framework to sequentially update the probability of change points. It excels in handling noisy signals by incorporating noise into the model, making it useful for online detection in streaming data \cite{adams2007}. However, BOCPD can be computationally intensive and sensitive to prior assumptions. Pruned Exact Linear Time (PELT), another widely-used method, minimizes a cost function plus a penalty term to avoid over-fitting to noise, offering scalability to large datasets and flexibility in cost function choice. Its performance depends on a penalty term, making it suitable for large datasets with unknown change points \cite{killick}.
For signals with high noise levels, Total Variation Regularization (TVR) effectively detects change points while simultaneously denoising, though it may smooth out small changes \cite{RUDIN1992259}. Cumulative Sum (CUSUM) detects mean shifts by monitoring the cumulative sum of deviations from a reference point, with robust variants for noisy data, though it is less effective for variance or distributional changes \cite{page}. Dynamic Programming (Exact Segmentation) offers an exact solution by minimizing the cost function over all possible segmentations but is computationally expensive and requires the number of switches to be known \textit{a priori}, making it more suitable for offline detection \cite{jackson}.
Hidden Markov Models (HMMs) are ideal for time series with hidden states and can model transitions as change points, though selecting the number of hidden states and computational cost are concerns \cite{rabiner}. Kernel-Based Methods, which map data to high-dimensional spaces, excel at detecting complex, non-linear changes but require parameter tuning and can be computationally intensive \cite{Harchaoui2007}. Lastly, Wavelet-Based Methods analyze signals at multiple scales, naturally denoising while detecting change points, making them suitable for noisy data with multiscale structures \cite{mallat}.
The choice of method depends on the type of change, noise level, and whether online or offline detection is needed. For low-to-moderate noise, PELT or BOCPD are recommended. For high noise, TVR or Wavelet-Based Methods work well, while for complex signals, Kernel Methods or HMMs are preferable.

In the field of switch detection, different cost functions are used to tailor the detection method to the specific characteristics of the data. $L_1$ and $L_2$ cost functions are basic measures used for minimizing absolute and squared residuals, respectively, where $L_1$ is more robust to outliers, while $L_2$ is sensitive to large deviations. The Normal cost from \cite{killick} assumes data follows a normal distribution and detects changes in mean and variance. More complex cost functions using Radial Basis Functions (RBFs) and sinusoids detect changes in data patterns by comparing similarity measures, with RBFs suitable for non-linear changes, and sinusoids focusing on angular shifts in data \cite{vert}.  Cost functions for data with linear relationships and continuous linear relationships account for gradual linear changes over time \cite{Harchaoui2007}.  Cost functions which rank data values are robust to non-Gaussian noise \cite{LIU201372}. Maximum likelihood cost functions adapt to various statistical models for flexible change detection \cite{Lavielle} . Finally, auto-regressive cost functions identify changes in time series data with auto correlation by capturing shifts in time dependencies \cite{davis06}. Every cost function mentioned can be selected based on the specific structure and noise characteristics of the data.

For the implementation of switch detection in our framework for the examples in this paper, we used binary segmentation (\texttt{binseg} which is a recursive algorithm that was devised to detect switches in noisy data \cite{Bai_1997}.  \texttt{Binseg} is sequential in nature, in that it first detects one change point, or switching point, in a given input data signal, then it splits the data at this change point, and this process is recrusively repeated on the two resulting sub-signals until no more change points are detected in all sub-signals \cite{wild_bin_seg}. 

Our implementation of \texttt{binseg} uses two hyperparameters where $SD_{algo}$ below is \texttt{binseg} in our examples.

\begin{align}
    \sigma &- \begin{cases}
    \text{standard deviation of noise,} & \text{if number of switches, $N_s$, is unknown,} \\
    \text{number of switches $N_s$,} & \text{if number of switches, $N_s$, is known.}
    \end{cases} \\
    s_g &- \text{switch gap : minimum number of data points between detected switches.} 
\end{align}

\begin{algorithm}[H]
\caption{Switch Detection Over All States of a System's Model}
\KwIn{$X_{data} \in \mathbb{R}^{m \times N}$: data matrix \\
$t \in \mathbb{R}^N$: independent variable \\
$s_g \in \mathbb{N}$: switch gap threshold \\
$SD_{cost}$: switch detection cost function model \\
$SD_{algo}$: switch detection algorithm \\
$\sigma > 0$ or $\sigma = N_s \in \mathbb{Z}_{\geq 0}$: noise level or fixed number of switches}
\KwOut{$N_s \in \mathbb{Z}_{\geq 0}$: number of switches detected \\
$\hat{t}_k$: switch locations (indices of $t$ where $k \in \{1,...,N_s\}$)}

\SetKwFunction{FMain}{Switch\_Detect}
\SetKwProg{Fn}{Function}{:}{}

\Fn{\FMain{$X_{data}$, $t$, $s_g$, $\sigma$, $SD_{cost}$, $SD_{algo}$}}{
    $(m, N) \gets \text{shape}(X_{data})$\;
    $\hat{t}_k \gets \emptyset$\tcp*{Initialize true switch locations across states}
    
    \For{$j \gets 1$ \KwTo $m$}{
        $\hat{t}_k \gets SD_{algo}(X_{data}[j, :], SD_{cost}, \sigma)$\tcp*{Detected and sorted switches for state $j$}

        \For{$i \gets 1$ \KwTo \text{length}$(\hat{t}_k)$}{
            \If{$\hat{t}_k[i] > \hat{t}_k[i-1] + s_g$}{
                $\text{result} \gets \text{append}(\hat{t}_k[i])$\tcp*{Applies $s_g$ for state j}
            }
        }

        \If{$\text{result}[\text{end}] == N$}{
            $\text{result} \gets \text{delete}(\text{result}, \text{end})$\tcp*{Last entry of $t$ cannot be a switch location}
        }

        $\hat{t}_k \gets \text{result}$\;

        \For{each $entry \in \hat{t}_k$}{
            \If{$entry \notin \hat{t}_k$}{
                $\hat{t}_k \gets \text{append}(entry)$ \tcp*{Removes redundant switches across states}
            }
            $\hat{t}_k \gets \text{sort}(\hat{t}_k)$\;
        }
    }

    \If{\text{length}$(\hat{t}_k) == 0$}{
        $N_s \gets 0$\;
        \Return $\hat{t}_k, N_s$ \tcp*{If no switches detected, return}
    }
 
    \If{\text{length}$(\hat{t}_k) > 0$}{ 
        \For{$i \gets 1$ \KwTo \text{length}$(\hat{t}_k)$}{
            \If{$\hat{t}_k[i] > \hat{t}_k[i-1] + s_g$}{
                $\text{result} \gets \text{append}(\hat{t}_k[i])$ \tcp*{Applies $s_g$ across states}
            }
        }
        $\hat{t}_k \gets \text{result}$\;
    }   
    
    \If{$\sigma == N_s \textbf{ and } \sigma \sim = length(\hat{t}_k)$}{
        throw error: ``detected $N_s \sim = $ expected $N_s$"
    }
    $N_s \gets \text{length}(\hat{t}_k)$\;

    \Return $\hat{t}_k, N_s$\;
}
\end{algorithm}\label{algo1}

Already existing as a feature of \texttt{binseg}, the user must input $\sigma > 0$ which is the standard deviation of the noise of the input signal. This noise level can be statistically estimated or inferred by the user as a pre-processing step.  Alternatively, the user may also iteratively try different $\sigma$ values until a desired model accuracy is achieved, or the user can smooth the data with a smoothing algorithm of one's choice such as the Savitzky-Golay filter \cite{Savitzky_Golay_1964}, FIR methods \cite{schnmid}, or Whittaker-Henderson smoother \cite{Whittaker_1922}.  The lower $\sigma$ is, the more sensitive the switch detection is, meaning that more switches will typically be detected for a given signal.  This is intuitive because \texttt{binseg} does not want to mistake noisy "jumps" for a switch in the signal.  There is also the case where the user may already know how many switches there are in the given input signal, and does not need \texttt{binseg} to estimate the number of switches.  If this is the case, the user may specify the number of switches to \texttt{binseg} as $\sigma = N_s$ and \texttt{binseg} will simply find the switch locations.  Note that $\sigma$ is the most important hyperparameter in determining the number of switches and their locations in given signals. 

We denote the "switch gap" as $s_g \in \N$ which does not allow there to be two parameter switches to be within $s_g$ discrete data points from each other.  This was implemented to allow some necessary user discretion when applying switch detection methods.  To expand, several methods such as \texttt{binseg}  may detect switches at consecutive time points and at other intervals that are too small for multiple switches to have occurred for certain classes of problems.  Simply setting $s_g=1$ allows for consecutive data points be detected as switches.   

For our implementation of \texttt{binseg}, all states (observed input signals) of an ODE or PDE are searched individually for the correct number of parameter switches as opposed to searching over all of the states at once.  We do this by running \texttt{binseg} on each state's signal, collecting every switch that occurs, and removing switches that have already been counted.  For example, suppose we have a two states of a system given by $x_1$ and $x_2$.  If it is determined that state $x_1$ experienced 2 switches at $t_1$ and $t_2$ and state $x_2$ experienced the same 2 switches, but also 2 more completely different switches at $t_3$ and $t_4$, then we say the system has 4 switches in total, namely 2 shared (coupled) switches, and 2 unique (uncoupled) switches in the second state (For an example of this, see Equation (\ref{eq:pvls2}).  This implementation grants the correct number of switches with accuracy when used with optimal switch detection algorithms on systems with discretely switching parameters.

%\begin{proposition}\label{prop: 2}
%     The implementation of switch detection in Algorithm 1 detects the correct number of parameter switches, $N_s$, with accuracy ($E_{N_s}=0$ \& $H_s=0$) for a dynamical system $\frac{dX}{dt} = f(X, p(t))$ where $X \in \R^{m}$ if the $SD_{algo}$ selected grants accuracy for each state in $X$.
%\end{proposition}
%
%\begin{proof}
%    Another way to phrase the proposition is that a system's parameter switches, $t_k$, will be equal to those detected, $\hat{t}_k$, if the specified $SD_{algo}$ detects all parameter switches across all states in $X$.  Let us then assume that the \textit{a priori} selected $SD_{algo}$ grants accuracy for each state in $X$.  In other words, all parameter switches that occur in a given state are accurately detected and accounted for.  Note that $t_k$ is a set of unique indices in which particular snapshots of $X_{data}$ indicate a switching parameter, and note that parameters can switch independently from state-to-state. If any of these parameter switches occur at the same locations across states, then our method removes these redundant parameter switches from the set of $\hat{t}_k$ so that each parameter switch is unique. So, if each state contains a set of parameter switches, then our algorithm returns of the union of the sets respectively obtained from each state.  Thus, we have that our method exhaustively accumulates all unique switches across all states with accuracy.  Hence,  we have shown that $\hat{t}_k = t_k$ .  
%\end{proof}

\begin{remark}
    If one was to naively use a switch detection algorithm directly on $X_{data}$ without the nuances of our implementation, then there would be severe errors $E_{N_s}$ \& $H_s$ in several cases such as the example given in Section (\ref{sec: non uni}) where the parameters switch at non-uniform intervals.  Ultimately, since systems with multiple states can have parameters that switch distinctly and independently across those states, the necessity of our dynamical systems switch detection implementation is clear.  
\end{remark}

\subsection{Optimization schemes, known static parameters, parameter bounds, constraints, and complexity}
\label{sec: schemes}

Upon experimenting with several numerical schemes to minimize our desired objective function.  Our best results came from robust, non-gradient based, deterministic optimization schemes, namely, the Nelder-Mead method \cite{nelder_mead} and the conjugate-based Powell's method \cite{powells_method}.  Auto-differentiation based methods and stochastic methods could often give comparable results, but not as consistently as the deterministic methods aforementioned.  The lack of consistency using auto-differentiation-based methods is most likely due to the combination of the general non-convex nature of objective functions and the stochastic elements of the methods.  It may be possible to get improved performance with auto-differentiation-based methods or in general from adding regularization to our objective function, but this is yet to be explored.  In \cite{Gabor} global optimization schemes are strongly recommended for handling the inherent non-convexity of these data-driven problems, especially in the context of biological systems.  With this insight, it is important to note that the modular nature of our framework allows for the implementation of global optimization schemes on finite domains instead of the optimization schemes mentioned previously.  Commonly used global optimization methods include basin hopping \cite{basin_hop}, brute force, differential evolution \cite{diff_evo}, dual annealing \cite{dual_anneal}, SHGO \cite{shgo}, and the DIRECT algorithm \cite{direct}. 
In the context of parameter estimation and especially when using global optimization methods, the utilization of parameter bounds and constraints is crucial to mitigate the computational complexity of navigating these spaces of parameters.
Parameter bounds and constraints can be used within our framework along with the distinctive feature of estimating a subset of parameters as constant.  
These features serve as ways to help tackle the non-convexity and computational efficiency of these data-driven optimization problems as they reduce the space of possible parameters when using global or local optimization schemes. 

 The computational complexity of an algorithm in our framework can be described as the computational complexity of the switch detection algorithm times the number of system states added to the estimated number of parameter switches times the complexity of the integration scheme ($\mathcal{O}(N)$) times the complexity of the optimzer.  This can be written with the  expression
\begin{align}
	m* \mathcal{O}(f_{SD}(N)) + \hat{N}_s * \mathcal{O}(N*f_{opt}(N, n_p)) \quad . 
\end{align}
where ($*$) represents standard multiplication, $f_{SD}$ is the complexity of the switch detection algorithm, and $f_{opt}$ is the complexity of the optimizer which depends on the number of data points and number of parameters.
Let us note that since the number of data points is known in the problem statement, the step-size between datapoints cannot change, thus making the complexity of the integration scheme $\mathcal{O}(N)$ .
Let us also note that in the most computationally complex scenarios that still achieve an optimal solution (convergence), we would have $\hat{N}_s \approx N$, \quad $f_{SD}(N) = N^2$ , \quad $ m \geq N $ .  From here we can see that the rightmost term typically dominates in most use cases.  In the examples we present the overall complexity is $\mathcal{O}(N^2)$ which is granted by the dominance of the rightmost term.  This also is the minimum complexity that can be achieved by algorithms within our framework.
  
\subsection{Sparse (LASSO) Regression for continuously varying parameters}
\label{sec: sparse id}

Sparse regression and dictionary methods have been successfully used across several domains to accurately fit parametric curves to signals in the context of dynamical systems \cite{brunton2016sindy}. Sparse signal reconstruction on fixed and adaptive supervised dictionary learning for transient stability assessment is provided by Dabou et al \cite{sparse_sig}. Devito et al. provide a dictionary optimization method for reconstruction of compressed ECG signals \cite{Devito2021}.  A framework for signal processing using dictionaries, atoms, and deep learning is given by Zhang and Van Der Baan \cite{signal_dict}.  Given the advances in using dictionaries for sparse regression, we use it for the purposes of identifying the structures of continuously changing parameters in parameter-varying systems. 

Suppose that we wish to capture a continuously changing parameter $p(t) : \R \mapsto \R $ as opposed to a piecewise-constant parameter that changes a discrete number of times $p = \{p_1, ... , p_n \} $ where $n$ is the number of segmented data intervals.
We can estimate $p(t)$ using our framework as constructed with some specifications on switch detection and an additional step of sparse regression to identify the continuously changing parameters.  To do this, we segment the data into $\frac{N}{6}$ intervals where $N$ is the number of data points in the input signal states.  This ensures that the parameters of a given system are not over-sampled or under-sampled.  Once our parameters have been sampled, we use the sparse LASSO regression to fit a function to each parameter where each function is a linear combination of a dictionary of functions.  Let us denote this dictionary as $\mathcal{D}_\theta (t)$, where each column of this matrix is a parameterized function denoted as $f_i(\theta_i,t)$, where $i \in \{ 1,...,r \}$, and $r$ is the number of dictionary functions, and each $\theta_i$ is a set of parameters for $f_i$.  The linear combinations of these dictionary functions shown approximate each parameter, $p(t)$. 

However, in practice, we cannot assume a well chosen dictionary, convexity of the objective function, uniqueness of a solution, or reasonable dimensionality. While these outlined challaenged are largely unresolved in the literature, we can use LASSO to promote solutions $\mathcal{D}_{\theta} (t) w$ that approximate $p(t)$ with minimal terms.  In an applied computational setting, let us denote $p$ as a finite vector approximation of the continuous time-varying parameter $p(t)$.  Then in order to solve 
\begin{equation}
    \min_{\theta, w} \| p - \mathcal{D}_{\theta} w \|_2 + \lambda \|w \|_0
    \label{eq:sparse}
\end{equation} 
where $\lambda \in \R$ is the sparse regularization parameter, and $w \in \R^r$ is the vector of weights.  The objective function seen in equation (\ref{eq:sparse}) allows us to achieve an approximation of the parameter $p \approx D_{\theta}w$ with only a sparse subset of our dictionary, granted to us by the $\lambda \| w \|_0$ term.  Note that that 1-norm for the term $\lambda \| w \|_1$ may grant better results when using gradient-based optimization methods and in other cases. 

\begin{figure}[H]
    \centering
    \includegraphics[width=.81\linewidth]{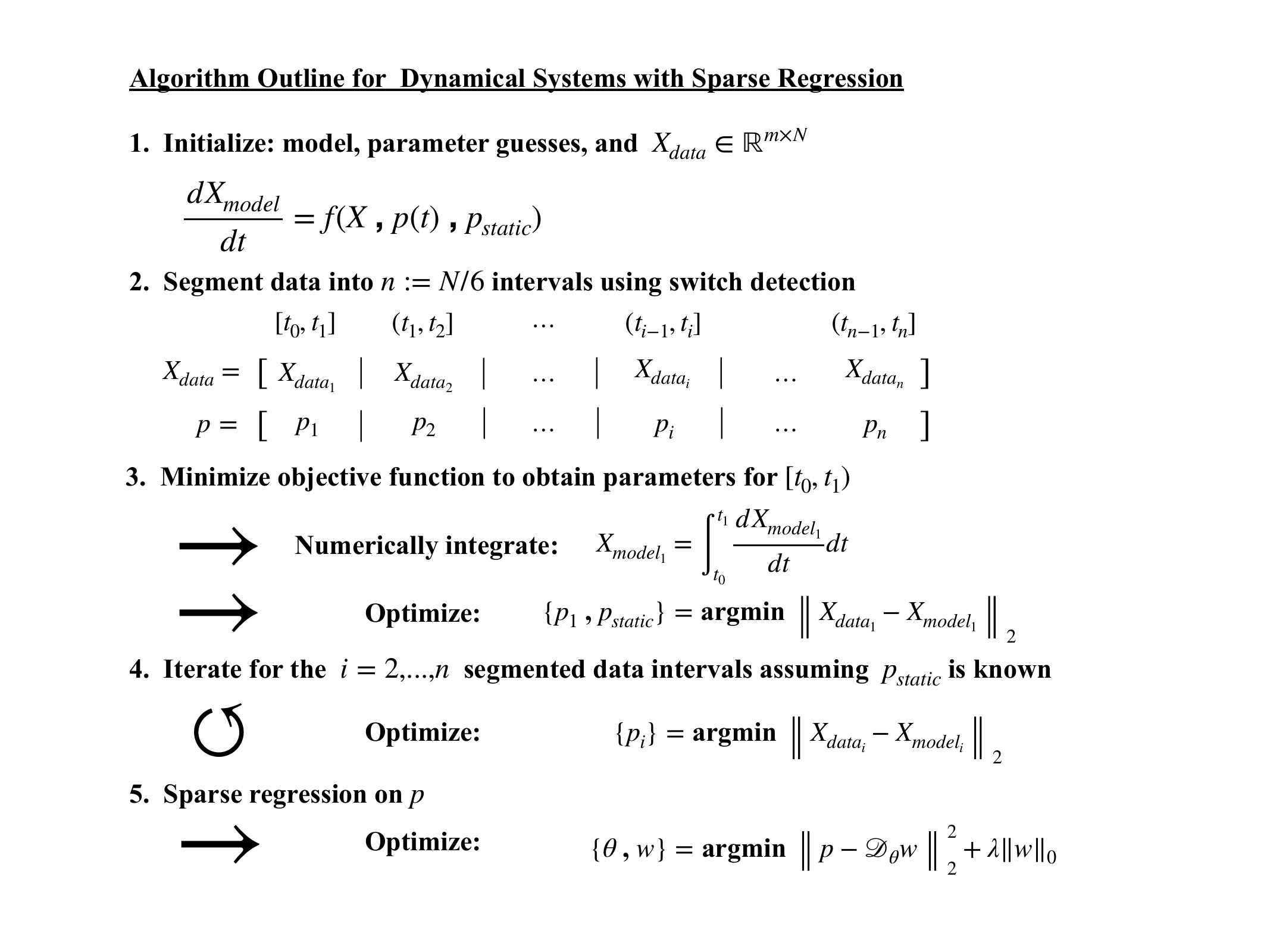}
    \caption{Outline of algorithm with sparse regression for the identification of continuously varying parameters.}
    \label{fig: algo pc}
\end{figure}

\section{Results and performance}

We first analyze our framework's performance on a nonlinear, parameter-varying system which has served as inspiration for the development of our algorithm. The model that for the Parameter-Varying Toggle Switch (PVTS) is given below and has been analyzed in \cite{Harrison}.  The parameters in this model have been shown to be globally unidentifiable in \cite{villaverde2016structural}, so our goal is to achieve accurate parameter estimates in spite of this. 
\begin{align}\label{eq:PVTS} 
    \begin{pmatrix} 
        \dot{x} \\  \dot{y}
    \end{pmatrix} &= \begin{pmatrix} 
        \frac{\alpha_1(t) }{1+(y/k_y)^{n_y} } -\delta_1(t) x \\[8pt]
        \frac{\alpha_2(t)}{1+(x/k_x)^{n_x}} -\delta_2(t) y 
    \end{pmatrix}  \\
    \alpha_1 = \alpha_2  & = 
    \begin{cases}
     1.0  &  \text{if } t \in [0, 12] \\
     8.0  &  \text{if } t \in (12, 24] \\
    \end{cases} \\
    \delta_1 = \delta_2  & = 
    \begin{cases}
     0.3  &  \text{if } t \in [0, 12]  \\
     0.6  &  \text{if } t \in (12, 24]  \\
    \end{cases} \\
    k_x = k_y  & = 1.0 \\
    n_x = n_y  & = 3.35 
\end{align}
where $x$ and $y$ represent the concentrations of repressor one and repressor two respectively; $\alpha_1(t)$ and $\alpha_2(t)$ denote the effective rates of synthesis of repressor one and two. $\delta_1(t)$ and $\delta_2(t) $ are the self decay rates of repressor one, repressor two; $n_x$ and $n_y$ represent the respective cooperativities of repression for promoter one and promoter 2; the Michaelis constants $k_x$ and $k_y$ are respective binding affinities.  

We simulated $N=1000$ data points from this model using LSODA as the ODE solver. The algorithm sub-steps we used were: \texttt{binseg} for switch detection with auto-regressive cost, $\sigma = 10^{-5}$, $s_g = 5$, Differential Evolution for a global parameter search where each parameter is in $[0,10.00]$ over a given time interval, Forward Euler for numerical integration of the model during the optimization step.  Our switch detection method achieved accuracy in terms of $H_s=0$ and $E_{N_s}=0$. When assuming $n_{x,y}$ and $k_{x,y}$ to be static parameters as described in Section (\ref{sec: schemes}), we have $E_p = 0.08$ and $E_t = 0.0041$ for the given PVTS as opposed to leaving this assumption out of our algorithm for a worse performance of $E_p = 0.2$ and $E_t = 0.0041$.  Thus, instantiating a subset of parameters as constant granted us more accurate parameter estimates than if we had assumed all parameters to be varying. 

\begin{figure}[!htb]
    \centering
    \includegraphics[width=.9\linewidth]{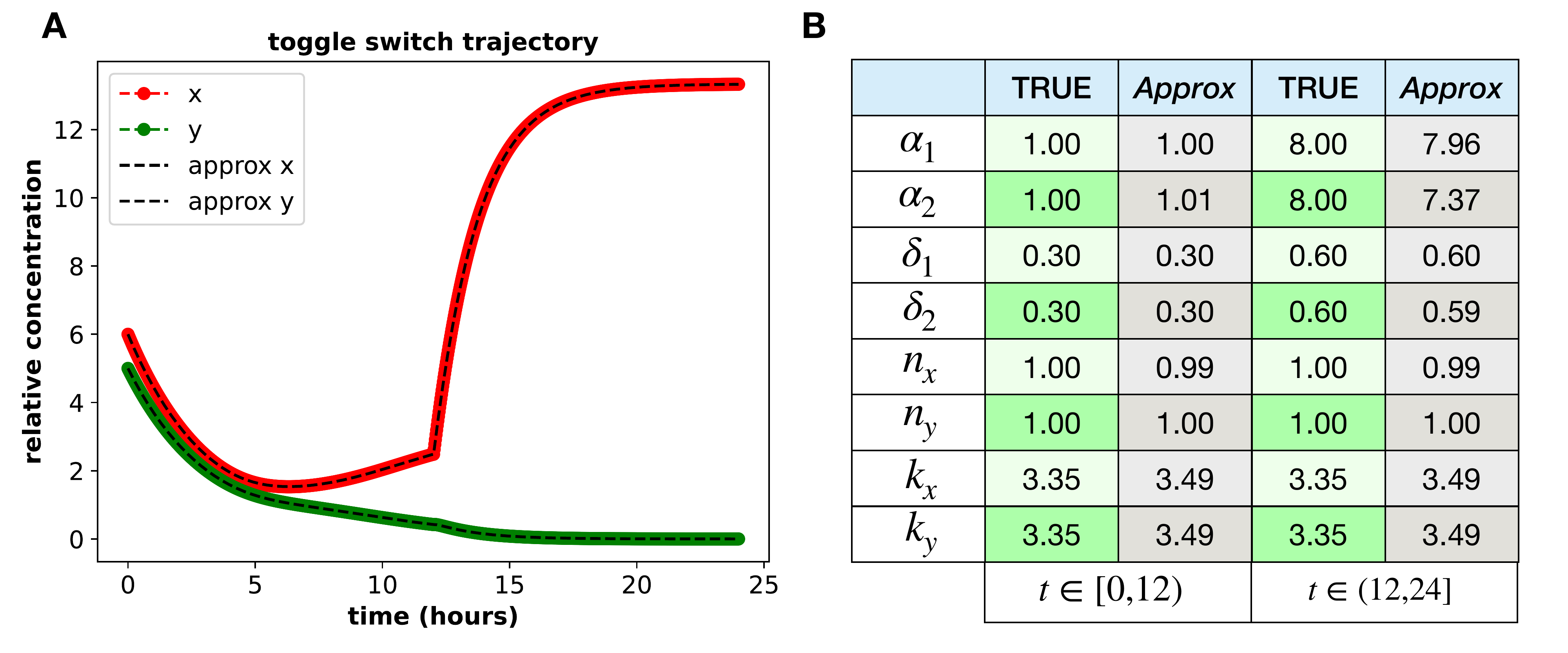}
    \caption{(A) Solution to parameter varying toggle switch and reconstruction of solution with estimated parameters. (B) Estimated parameters for parameter varying toggle switch. }
    \label{fig: toggle switch}
\end{figure}

We see that we are able to accurately estimate the static and bi-phasic parameters of the PVTS as desired.  Achieving this in our framework was non-trivial in that the parameter space contains several local optima over each time interval.  To get obtain the level of accuracy we have shown, we employed the global optimization scheme Differential Evolution to minimize $E_t$, and we made the \textit{a priori} instantiation of the $n_{xy}$ and $k_{xy}$ parameters to be estimated as constant functions.  Using Nelder-Mead or Powell's methods as opposed to the Differential Evolution gave us worse errors of $E_t = 0.0049$ and $E_p = 0.79$.       

\subsection{PDEs and parametric curves}
\label{sec: pdes and pcs}

To illustrate how our framework can handle PDE's, we set up the pertinent algorithmic sub-steps and we ran the algorithm on the heat equation  where the diffusion coefficient is the parameter varying function $D(t) = p_1 (t) t^{p_2(t)}$ as shown in Equation (\ref{eq: pde test}) where $p_1(t)$ and $p_2(t)$ are piece-wise functions in time.  A physical understanding of this phenomenon is provided in Appendix~\ref{App1}.
  
\begin{align}\label{eq: pde test} 
    \quad \frac{\partial}{\partial{t}} u(t,x) & =  p_1(t) t^{p_2(t)} \frac{\partial^2}{\partial{x^2}} u(t,x) \\
    \text{where \quad }    p_1(t)  & = 
    \begin{cases}
     0.1  &  \text{if \quad } t \in [0, 0.25] \\
     2.1  &  \text{if \quad} t \in (0.25, 0.5] \\
     4.1  &  \text{if \quad} t \in (0.5, 0.75] \\
     6.1  &  \text{if \quad} t \in (0.75, 1] \\
    \end{cases} \\
    p_2(t)  & = 
    \begin{cases}
     0.5  &  \text{if \quad} t \in [0, 0.25] \\
     0.4  &  \text{if \quad} t \in (0.25, 0.5] \\
     0.3  &  \text{if \quad} t \in (0.5, 0.75] \\
     0.2  &  \text{if \quad} t \in (0.75, 1] \\
    \end{cases} \\
    \text{with initial condition \quad } u(t=0, x) & = \begin{cases}
     x  &  \text{if \quad} x \in [0, \frac{\pi}{2}] \\
     \pi - x  &  \text{if \quad} x \in (\frac{\pi}{2}, \pi] \\
    \end{cases}\label{eq: ic heat} \\
    \text{and with boundary condition \quad } u(t, x=0) & = u(t, x=\pi) =0 \quad . 
\end{align}

\begin{figure}[H]
    \centering
    \includegraphics[width=.9\linewidth]{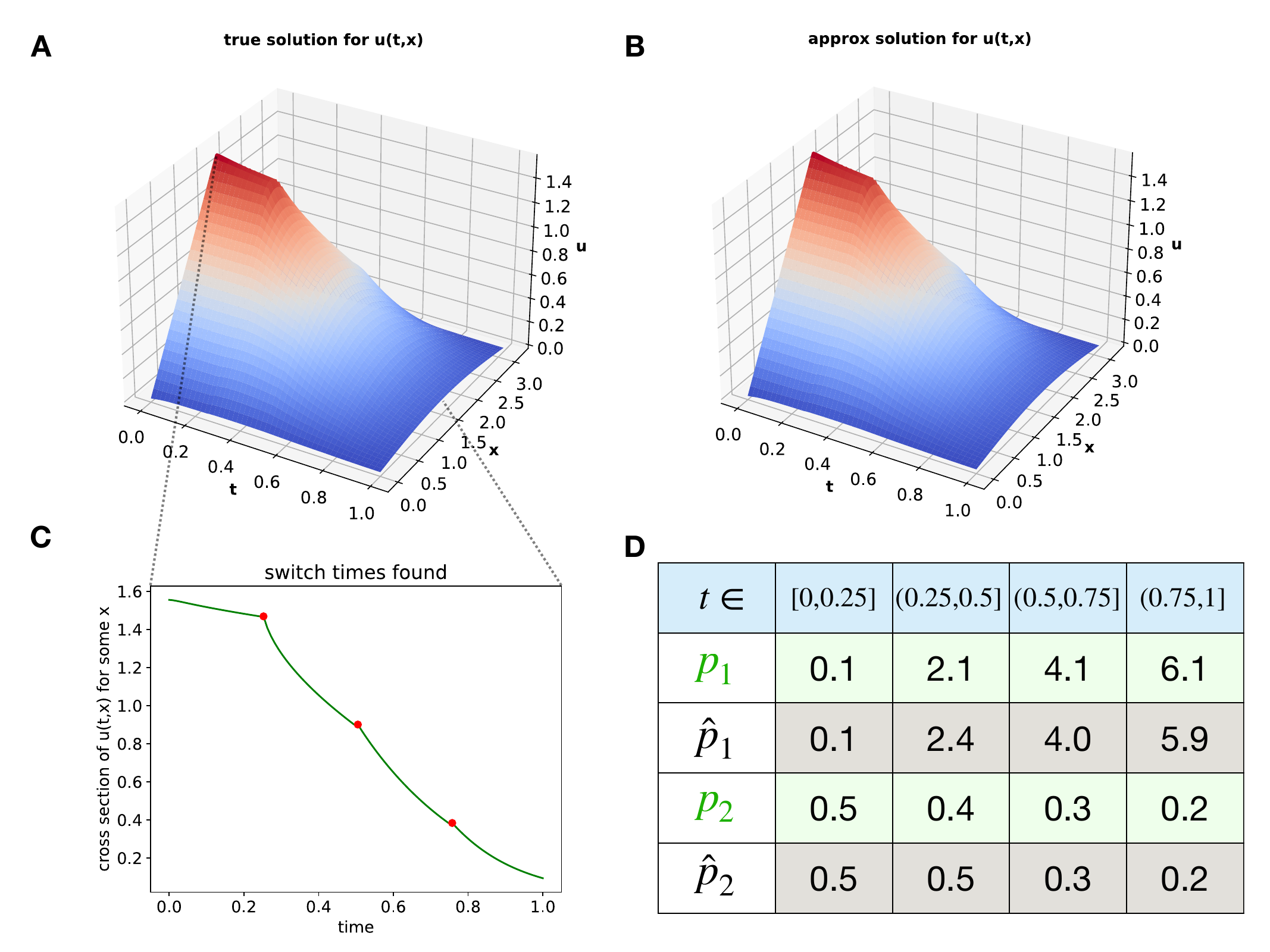}
    \caption{(A) True solution to heat equation PDE. (B) Reconstruction of solution to heat equation PDE with estimated parameters . (C) Switches detected from cross-sectional slice of PDE data. (D) Estimated parameters for heat equation PDE from data.}
    \label{fig: pde}
\end{figure}

We used the method of lines (discretization of spatial derivatives) as explicitly written in Equation (\ref{eq: pde test}) to represent the PDE as a system of ODEs with 100 uniformly spaced grid points in $x$ $(M=100)$ and 100 time points over the domain $(t,x) = [0,1] \times [0,\pi]$.  The ODE solver in the optimization sub-step was Forward Euler, we used \texttt{binseg} for switch detection with auto-regressive cost, $\sigma = 10^{-3}$, $s_g = 10$, and Nelder-Mead for locally optimal parameters with no pre-set bounds.
We see that we have achieved $H_s = 0$, $E_{N_s} = 0$, $E_p = 0.19$, and $E_t = 0.001$ via the algorithm specifications.  Thus, we have shown that that an algorithm specified within our framework can accurately estimate the varying parameters of a PDE.

\subsection{Robustness of Switch detection and parameter estimation given noisy data}

Quantifying the noise-robustness of our framework is needed given the ever-increasing need to process and understand noisy data for system identification and control \cite{sakib2024pi, sakib2024hankel, sinha2023online}.  To quantify the robustness of our algorithm, we added white noise (zero mean and finite variance) to synthetic data that was generated from three different known ground-truth models which all have distinct orders, numbers of parameters, and system structures.  We then ran this noisy data through our algorithm and identified relationships between the noisiness of the data and the types of error described in Table (\ref{tab: error types}) which are inherent to our algorithmic substeps.

To illustrate, one of the ground truth model we used is a biological system given in Equation (\ref{eq:pvls}) where $m$ is mRNA concentration, $p$ is protein concentration, $\alpha_{m,p}$ are respective rates of synthesis, and $\delta_{m,p}$ are the respective degradation rates.  This is a standard protein synthesis model in which mRNA is transcribed into a given protein.  Note that in this model, we have that the protein production is dependent on the transcription of mRNA, but not vice versa.  Hence, the model is in agreement with the central dogma of biology.  

A motivating factor for this work was to validate novel genetic circuit models with time-varying kinetic rates such as the one below.  Based on the equations provided, it is clear that the $\alpha_{m,p}$ rates of synthesis are the varying parameters of the model, and that the respective degradation rates are modeled as unchanging.  In an experimental setting, the data that would be collected would be subjected to some measurement noise from the instruments that are used to quantify time-series gene expression levels. 

% take out gray box from this figure
\begin{figure}[H]
    \centering
    \includegraphics[width=.8\linewidth]{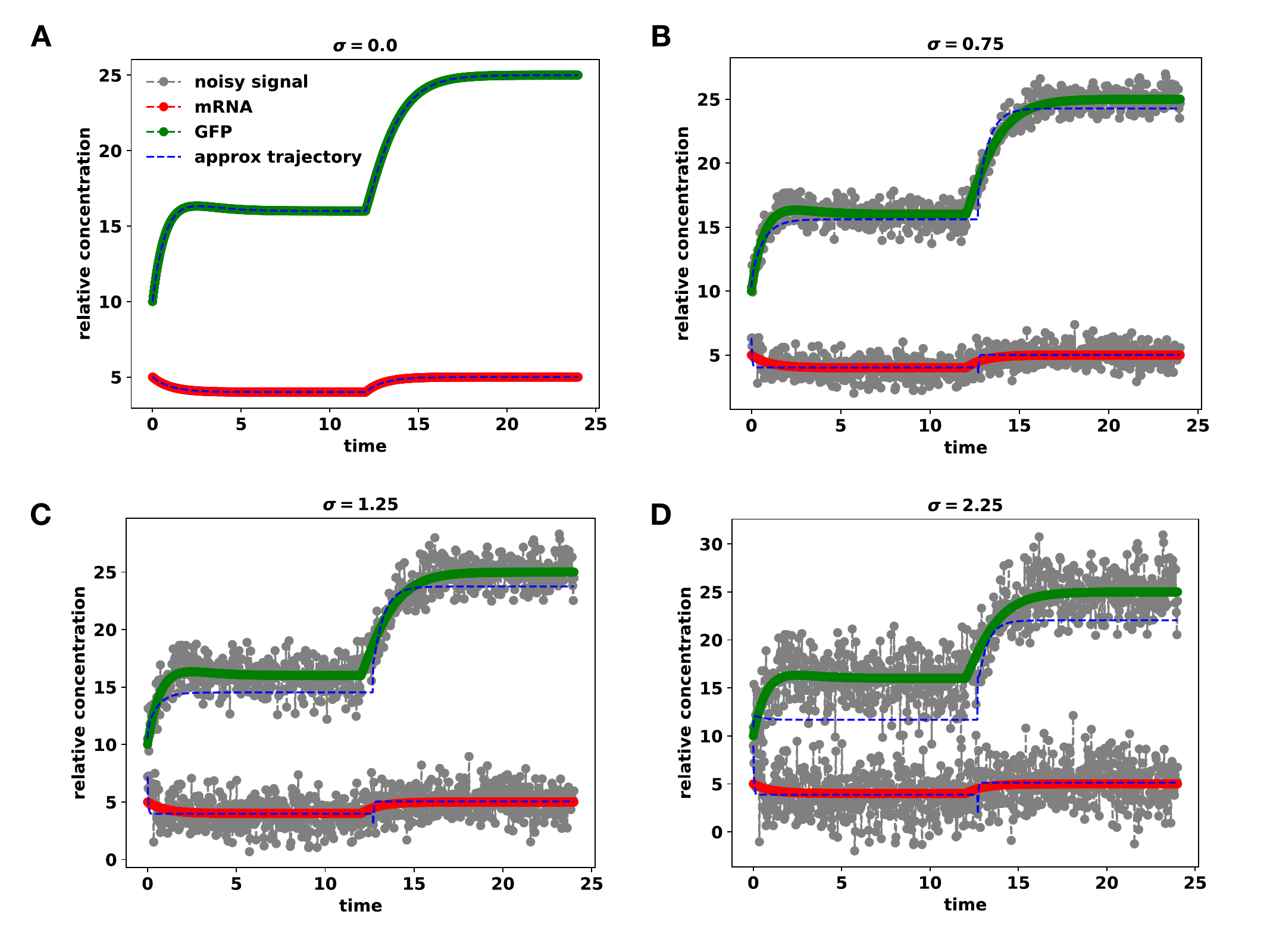}
    \caption{(A) Synthetic data and reconstruction from estimated parameters. (B) Synthetic data with white noise $(\sigma =0.75) $ and reconstruction from estimated parameters. (C) Synthetic data with white noise $(\sigma =1.25) $ and reconstruction from estimated parameters (D) Synthetic data with white noise $(\sigma = 2.25) $ and reconstruction from estimated parameters. }
    \label{fig: noise}
\end{figure}

\begin{align}\label{eq:pvls}
\dot{m} &= \alpha_m (t) - \delta_m (t) m \\
\dot{p} &= \alpha_p (t) m - \delta_p (t) p  \\
\alpha_m(t) & = \alpha_p(t)  = 
    \begin{cases}
     4.0  &  \text{if \quad} t \in [0, 12] \\
     5.0  &  \text{if \quad} t \in (12, 24] 
    \end{cases} \\
    \delta_m(t) & = \delta_p(t)  = 
     1 \quad \forall t \in [0, 24]  \quad .
\end{align}

Let the observables of our synthetic data be denoted as $x := [m \quad p]^T$ and let the white noise we add to it be denoted by $v_w$.  Each element of the white noise vector $v_w$ is sampled from a normal distribution with zero mean and finite standard deviation, $\sigma$.  Our algorithm was used to estimate the time-varying parameters of the given model (\ref{eq:pvls}) using $X_{data}:=x + v_w$.  A subset of the these runs are found in Figure (\ref{fig: noise}), where we manually set the number of switches to be detected in the data to be $N_s =1$ so that we could measure noise effects in this idealized best case scenario.  In the worst case scenario where the number of parameter switches, $N_s$, is unknown, switch detection algorithms often mistake ``jumps" in the data from the noise as switches which leads to worse performance than what we show here.  This worst case scenario often gives relatively accurate data reconstruction, but unreliable parameters and untrustworthy detection of switches.

The ODE solver we used for the numerical integration in the optimization sub-step was Forward Euler and the algorithm sub-steps we used were: \texttt{binseg} for switch detection with auto-regressive cost, $\sigma = N_s = 1$, $s_g = 10$, Nelder-Mead for locally optimal parameters with no pre-set bounds.  We used Forward Euler here because LSODA performed worse than Euler’s method within computing the objective function.  Time integration errors are likely to affect the estimation, since the loss function depends on the rolled-out solution. To reduce such errors, one can utlize the modular nature of our framework to implement other ODE solvers just as we have done.

\begin{figure}[!htb]
    \centering
    \includegraphics[width=\linewidth]{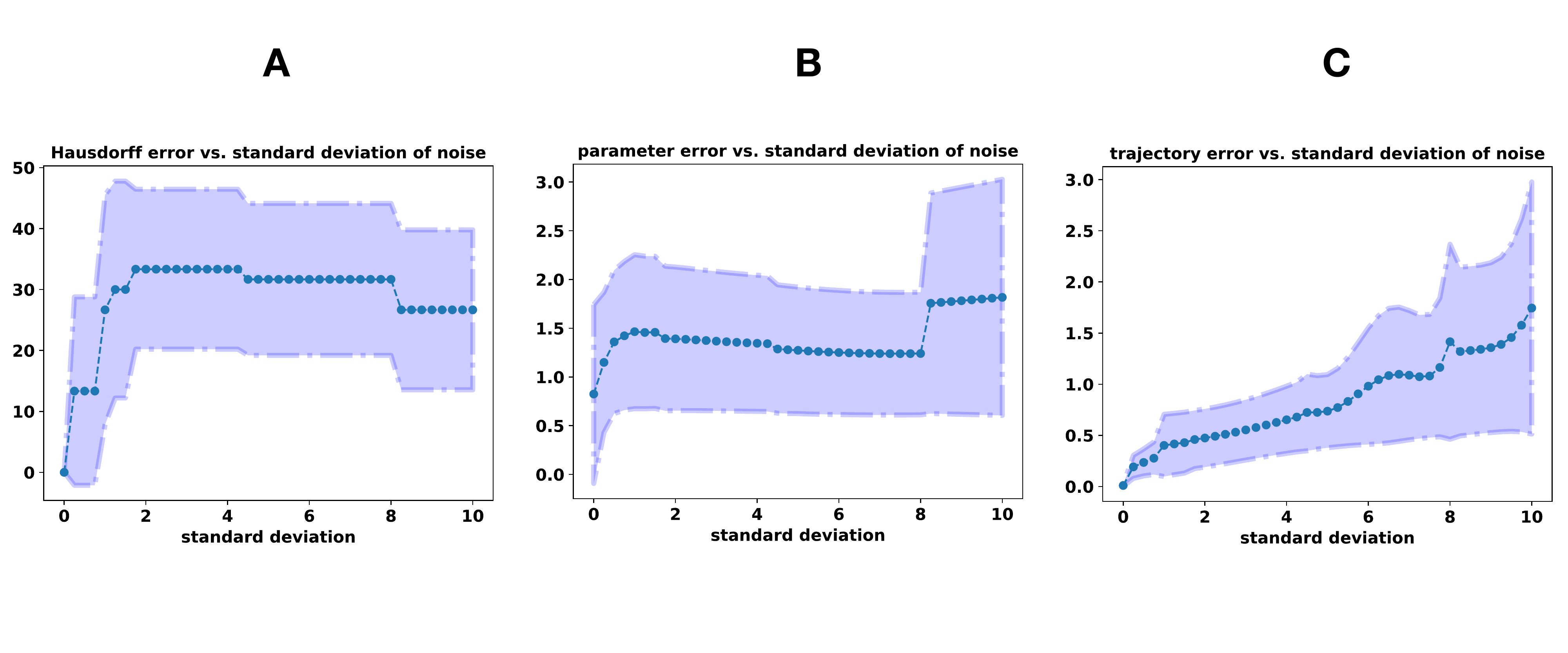}
    \caption{(A) Hausdorff error vs. standard deviation of noise.  (B) Parameter error vs. standard deviation of noise. (C) Trajectory error vs. standard deviation of noise. Dotted lines are the averages over several runs and dynamical systems.  Shaded regions are standard deviation error bars computed from all runs}
    \label{fig: noise error}
\end{figure}

From both Figure (\ref{fig: noise}) and Figure (\ref{fig: noise error}), it is clear that our method for parameter estimation shows increases in all error metrics with with increases to the noise in data. There is a cascading effect of error propagation due to noise $(v_w)$ : \texttt{binseg} inaccurately detects switches which leads to inaccurate data segmentation $ \left (\uparrow D_s \text{ and } \uparrow E_{N_s} \right )$ which causes inaccurate parameter fitting $(E_p)$, and thusly, parameterized model trajectories inaccurately approximate data $(E_t)$.  This cascading effect is summarized as

\begin{align}\label{eq: err prop}
   \uparrow \| v_w \| \implies \left ( \uparrow H_s \text{ and } \uparrow E_{N_s} \right ) \implies \uparrow E_p \implies \uparrow E_t \quad . 
\end{align}

This cascading propagation of error is one of the drawbacks of algorithms in our framework.  Methodological sub-steps specializing in the robust handling of noisy data can be added on at various points within our framework for improvements in this weak-point, but this is outside the scope of this paper and may be further expanded on in future work. 

\subsection{Switching parameters at non-uniform intervals}
\label{sec: non uni}

The special handling of our adapted switch detection framework to detect discretely switching parameters in dynamical systems is heavily motivated by cases in which parameters across the states of a dynamical systems are switching independently from each other.  
Suppose now that we have the same system as in Equation (\ref{eq:pvls}) but with parameters that switch non-uniformly on different time intervals as depicted in Equation (\ref{eq:pvls2}). 

\begin{align} \label{eq:pvls2}
\dot{m} &= \alpha_m (t) - \delta_m (t) m \\
\dot{p} &= \alpha_p (t) m - \delta_p (t) p  \\
\text{where \quad }    \alpha_m(t)  & = 
    \begin{cases}
     4.0  &  \text{if \quad} t \in [0, 5] \\
     5.0  &  \text{if \quad} t \in (5, 15] \\
     6.0  &  \text{if \quad} t \in (15, 24] \\
    \end{cases} \\
    \alpha_p (t)  & = 
    \begin{cases}
     4.0  &  \text{if \quad} t \in [0, 5] \\
     5.0  &  \text{if \quad} t \in (5, 10] \\
     6.0  &  \text{if \quad} t \in (10, 15] \\
     7.0  &  \text{if \quad} t \in (15, 20] \\
     8.0  &  \text{if \quad} t \in (20, 25] \\
    \end{cases} \\
    d_m(t) = d_p(t)  & = 
    \begin{cases}
     1  &  \text{if \quad } t \in [0, 25]  \quad .
    \end{cases}
\end{align}

\begin{figure}[H]
    \centering
    \includegraphics[width=\linewidth]{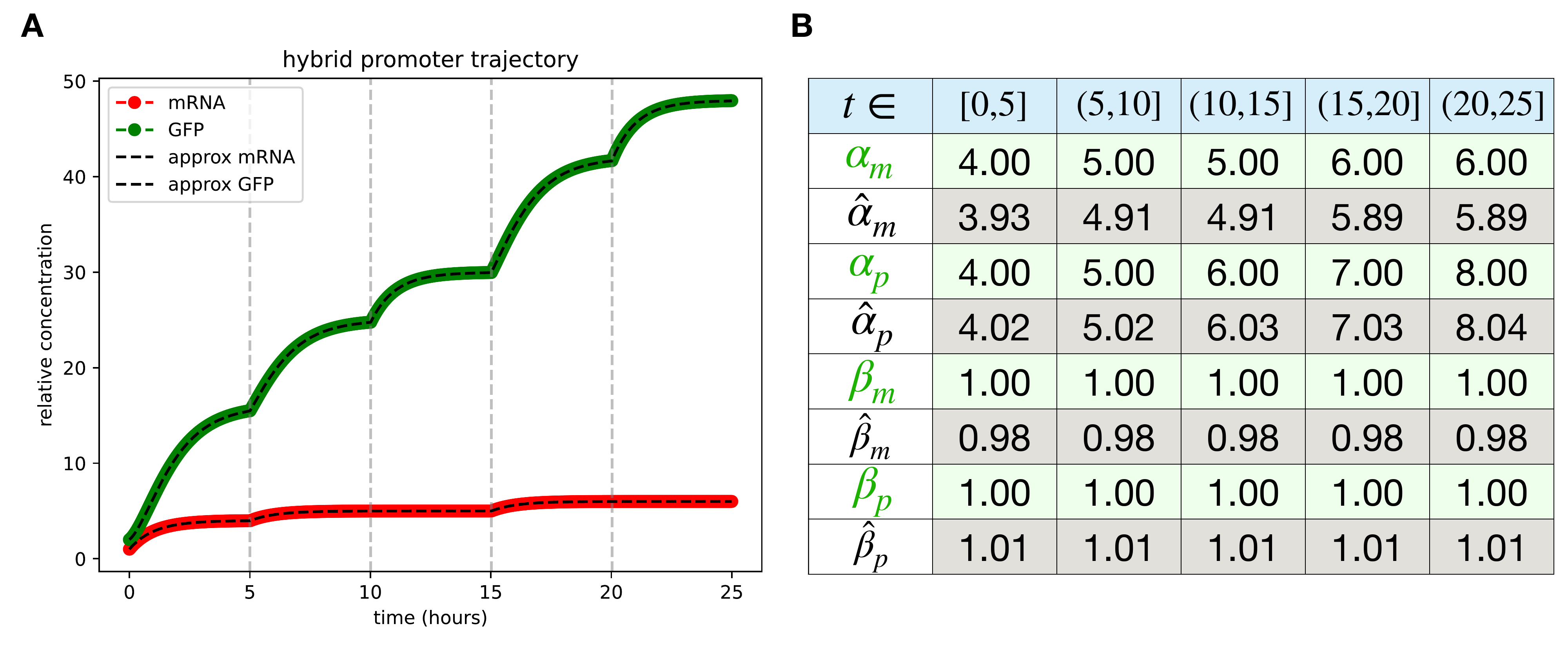}
    \caption{(A) Trajectory of hybrid promoter model with different parameters ($\alpha_m(t)$ and $\alpha_p(t)$) switching at different time intervals.  Switch locations are shown as grey dotted lines (estimated=true) (B) Approximated parameters of hybrid promoter model in which parameters switch at different time intervals.}
    \label{fig: diff time ints}
\end{figure}

We test our algorithm on System (\ref{eq:pvls2}) to show that our adaptation of \texttt{binseg} is able to handle these scenarios as shown in Figure (\ref{fig: diff time ints}).  The ODE solver we used for the numerical integration in the optimization sub-step was Forward Euler and the algorithm sub-steps we used were: \texttt{binseg} for switch detection with auto-regressive cost, $\sigma = N_s = 1$, $s_g = 10$, Nelder-Mead for locally optimal parameters with no pre-set bounds.  With this algorithmic set-up, we achieved $H_s =0$, $E_{N_s}=0$, $E_t=0.001$ and $E_p = 0.05$.

\subsection{Sparse dictionary regression for continuously varying parameters}

\begin{figure}[H]
\centering
    \includegraphics[width=0.99\linewidth]{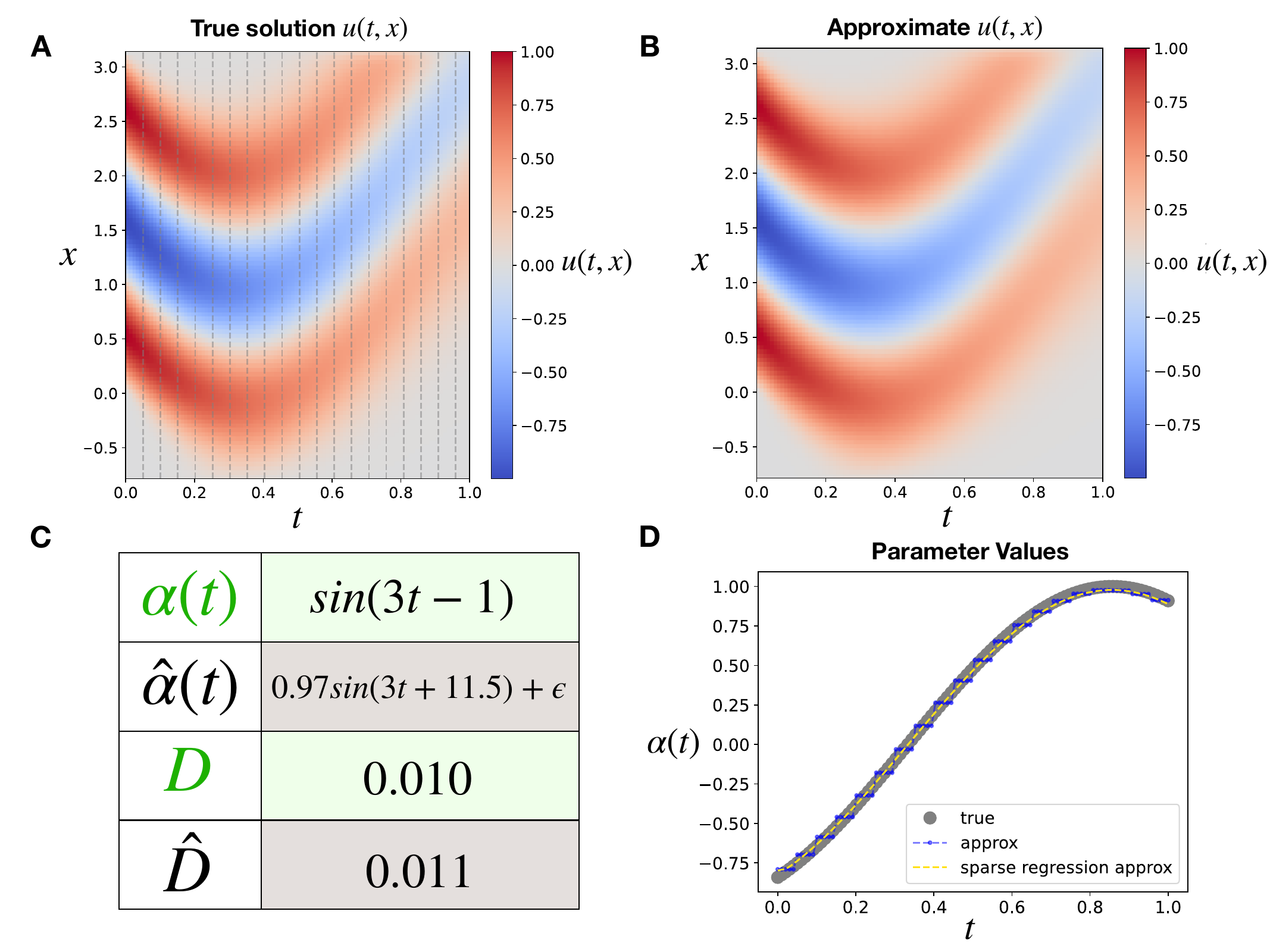}
    \caption{(A) Synthetic data from advection-diffusion equation containing a continuously varying parameter.  The vertical gray lines indicate the segmentation of the data for the sampling of the continuously varying parameter. (B) Reconstruction of heat equation data using estimated model parameters from algorithm with sparse regression. (C) Table of parametrized true and estimated parameters from the set of dictionary functions. (D) True and approximate parameters from our algorithm with and without sparse regression.}
    \label{fig:sparse}
\end{figure}

To test our method outlined in Section (\ref{sec: sparse id}), we estimate the time-varying parameter $\alpha(t)$ in the PDE model for advection-diffusion found in Equation (\ref{eq: sparse test}) using an algorithm specified within our framework adjoined with sparse regression. The main differences between this example and that of Section (\ref{sec: pdes and pcs}) are that we now have an advection term in the PDE denoted by the first order spatial derivative and that we are not assuming knowledge of the structure of our varying parameter,  $\alpha(t)$.  We proceed in the estimation of $\alpha(t)$ as a continuously varying parameter for two reasons: due to the horrendous accuracy obtained if estimating all parameters as constant, and the lack of coherent discrete switches detected in the data.  Intuitively, it is clear from the data found in Figure (\ref{fig:sparse}) that the direction and velocity of advection is changing.  We assumed in our test run that the diffusion coefficient $D$ was to be estimated as a constant parameter. 

\begin{align}\label{eq: sparse test} 
    \quad \frac{\partial}{\partial{t}} u(t,x) & = \alpha(t) \frac{\partial}{\partial{x}} u(t,x)  +  D \frac{\partial^2}{\partial{x^2}} u(t,x) \\
    \text{where \quad } \hspace{10mm} \alpha(t) &= sin(3t - 1) \text{\quad and \quad } D  = 0.01 \\
    \text{with initial condition \quad } u(0, x) & = \begin{cases}
     sin(3x)  &  \text{if \quad} x \in [0, \pi] \\
     0  &  \text{if \quad} \text{else} \\
    \end{cases} 
\end{align}

%\begin{equation} \label{eq: dict test}
%     \mathcal{D}_\theta (t) = \begin{bmatrix}
%        | & | & | \\ 
%        e^{\theta_1 t^{\theta_2}} & sin(\theta_3 t + \theta_4) & t^{\theta_5} \\ 
%        | & | & | \\ 
%    \end{bmatrix} \quad . 
%\end{equation} 

We used the method of lines to set up this PDE and we used LSODA to simulate the data with 100 uniformly spaced grid points in $x$ $(M=100)$ and 100 time points over the domain $(t,x) = [0,1] \times [-\pi/4,\pi]$.  The ODE solver we used for the numerical integration in the optimization sub-step was Forward Euler and the algorithm sub-steps we used were: \texttt{binseg} for data segmentation with auto-regressive cost, the number of assumed switches as $\sigma = N_s = N/6$, the switch gap as $s_g = 1$, Nelder-Mead for locally optimal parameters with pre-set bounds such that $|\alpha(t)| < 10$ and $|D| < 10$, and sparse regression with regularization term as $\lambda \| w \|_1$ where $\lambda = 0.01$.  With this algorithmic setup, we achieved $E_t=0.001$ and $E_p = 0.001$.

Without our implementation of sparse regression, we achieved a parameter error of $E_p = 0.16$.  With sparse regression we achieved a parameter error of $E_p = 0.05$, thus showcasing an improvement in parameter estimation accuracy through the utilization of sparse regression.  Further, sparse regression has allowed us to reconstruct the parameterization of $\alpha(t)$ as $\hat{\alpha}(t) = 0.97sin(3t + 11.5) 
 + \epsilon$ where $\epsilon$ contains the minuscule terms from the other dictionary functions, namely, $\epsilon = (10^{-13})e^{-2.55t^{0.6}} +0.01t^{0.96}$ .  Please also note that the offset of $-1$ radians in $\alpha(t)$ is approximately equal to the $11.5$ radian offset in $\hat{\alpha}(t)$.

\section{Discussion and Conclusions}

Currently, scenarios which contain only discrete parameter switches with accurate estimations of the switch times as well.  Scenarios which have only continuously varying parameters are able to be estimated accurately with the use sparse regression on the sampled parameter estimates as shown.  Working with parameters that are combinations of both discretely switching and continuously varying is successful up to the following extent: we can approximate these parameters/inputs reasonably well by sufficiently segmenting the data, but we cannot estimate the switch locations during this process.  In the broader context of the literature, our method is the first to attempt to estimate the swtich locations of parameters while also estimating the parameters themselves.  Upon reading our results, one might have thought to try running the dynamic programming switch detection algorithm  where $N_s$ is increased every iteration until the trajectory error is below some threshold. While this idea is somewhat sound upond reading, this method suffers greatly in practice with real data due to the following: (1) If $N_s$ is unknown then the user must have some insight into what error threshold to pick to give a reasonable estimate of $N_s$.  (2) This method is extremely computationally costly and often gives unreasonably high $N_s$ estimates with real, noisy data. (3) In general one is better off tweaking the hyper parameters of a switch detection algorithm like \texttt{binseg} in a practical setting with real, noisy data. 

In conclusion, we present a modular framework for estimating varying parameters in dynamical systems, focusing initially on the sub-class of problems where parameters switch discretely and are modeled as piecewise constant functions. By combining the strengths of switch detection, numerical integration, and optimization into a cohesive structure, our framework offers accuracy and flexibility to accommodate user-specific requirements. Our approach has been validated across diverse examples, including biological systems and physical models, demonstrating its adaptability to scenarios with fixed and switching parameters, non-uniform switching, and time-varying dynamics. Moreover, we characterized the framework’s robustness to measurement noise and extended its capabilities to continuously varying parameters using dictionary-based sparse regression with trigonometric and polynomial functions. This work establishes a foundation for addressing a wide range of parameter-varying system problems and provides a versatile tool set for advancing research in dynamical system modeling and control.

\section*{Acknowledgments}

For insightful discussions related to this work, the authors would like to thank Charles Johnson, Yanran Wang, Taishi Kotsuka, Harnoor Lali, SIAM, and Jo$\Tilde{\text{a}}$o Hespanha.  This work was funded in part by an NSF CAREER Award 2240176, the Army Young Investigator Program Award W911NF2010165 and the Institute of Collaborative Biotechnologies/Army Research Office grants W911NF19D0001,  W911NF22F0005, W911NF190026, and W911NF2320006.  This work was also supported in part by a subcontract awarded by the Pacific Northwest National Laboratory for the Secure Biosystems Design Science Focus Area “Persistence Control of Engineered Functions in Complex Soil Microbiomes" sponsored by the U.S. Department of Energy Office of Biological and Environmental Research.

\section*{Additional information}

We report no conflicts of interest. 

The code used to produce all data and figures in this paper can be found at the following \href{https://github.com/Jamiree/Estimating-Varying-Parameters-in-Dynamical-Systems}{Github repository}.  We also provide other supplementary examples not covered in this paper. 

\section*{Appendix A.1}\label{App1}
The theoretical physics underlying this example problem provides valuable insight into the behavior of complex materials under dynamic thermal conditions. In this case, the 1-D heat equation models the temperature distribution along a thin rod with an initial heat profile \( u(0,x) \) as defined in Equation (\ref{eq: ic heat}). As time progresses, the parameter representing the rod's heat diffusivity, \( D(t) \), evolves dynamically. This behavior reflects the rod's composition—a non-homogeneous composite material consisting of three distinct compounds, each with unique thermal properties and phase-change characteristics. 

As the rod cools, the heat diffusivity of one compound undergoes a significant change at its critical temperature, corresponding to the first parameter switch. Similar transitions occur for the other two compounds, resulting in the discrete parameter switches observed in the measured data. Assuming that previous studies have shown that the heat diffusivity of such composites changes exponentially with cooling, this behavior can be represented by the parameter structure \( D(t) = p_1(t) t^{p_2(t)} \). The non-homogeneous nature of the composite material justifies the discrete switching observed in the heat diffusivity.

\end{document}